\magnification =\magstep1
\overfullrule=0pt

\font\footfont=cmr8
\font\footitfont=cmti8

\font\large=cmbx12 at 14pt
\font\bb=msbm10

\def\cq{c(q)}
\def\xq{\delta(q)}
\def\tq{\Delta(q)}

\def\g{g}
\def\h{h}
\def\k{k}
\def\dg{\deg}

\def\lcm{\hbox{lcm}\hskip 0.3em}

\def\bbZ{\hbox{\bb Z}}
\def\Fq#1{{#1}^{[q]}}

\def\qed{\hfill\hbox{\vrule\vtop{\vbox{\hrule\kern1.7ex\hbox{\kern1ex}
  \hrule}}\vrule}\vskip1.5ex}
\def\eqed{\hbox{\quad
  \hbox{\vrule\vtop{\vbox{\hrule\kern1.7ex\hbox{\kern1ex}\hrule}}\vrule}}}
\def\mod{\hskip 0.3em\hbox{modulo}\hskip 0.3em}
\def\mvrule{\vrule height 3ex depth .4ex width .01em}

\def\margin#1#2{\rlap{\hbox to #2em{\ \hfill}\hskip-2em*#1*}}

\def\label#1{{\global\edef#1{\the\sectno.\the\thmno}}}
\def\pagelabel#1{{\global\edef#1{\the\pageno}}}

\def\isnameuse#1{\csname #1\endcsname}

\def\issecond#1#2{#2}
\def\isifundefined#1#2#3{
        \expandafter\ifx\csname #1\endcsname\relax #2
        \else #3 \fi}
\def\pageref#1{\isifundefined{is#1}
        {{\bf ??}\message{Reference `#1' on page [\number\count0] undefined}}
        {\edef\istempa{\isnameuse{is#1}}\expandafter\issecond\istempa\relax}}

\def\today{\ifcase\month\or January\or February\or March\or
  April\or May\or June\or July\or August\or September\or
  October\or November\or December\fi
  \space\number\day, \number\year}

\newcount\sectno \sectno=0
\newcount\thmno \thmno=0
\def \section#1{\vskip 0.3truecm\relax
        \global\advance\sectno by 1 \global\thmno=0
        \noindent{\bf \the\sectno. #1} \vskip 0.2truecm\relax}
\def \thmline#1{\vskip .2cm\relax
        \global\advance\thmno by 1
        \noindent{\bf #1\ \the\sectno.\the\thmno:}\ \ %
        \bgroup \advance\baselineskip by -1pt \it
        \abovedisplayskip =4pt
        \belowdisplayskip =3pt
        \parskip=0pt
        }
\def \wthmline#1{\vskip .2cm\relax
        \noindent{\bf #1:}\ \ %
        \bgroup \advance\baselineskip by -1pt \it
        \abovedisplayskip =4pt
        \belowdisplayskip =3pt
        \parskip=0pt
        }
\def \dthmline#1{\vskip 6pt\relax
        \noindent{\bf #1:}\ \ }

\def \thm{\thmline{Theorem}}

\def \endb{\egroup \vskip 0.2cm\relax}

\def \lemma{\thmline{Lemma}}

\def \conj{\wthmline{Conjecture}}

\def \note{\dthmline{Note}}
\def \prop{\thmline{Proposition}}
\def \proof{\vskip0.0cm\relax\noindent {\sl Proof:\ \ }}

\ %
\vskip 3ex
\centerline{\large Computations with Frobenius powers}

\vskip 4ex
\centerline{ Susan M. Hermiller
\unskip\footnote{*}{
\footfont The first author acknowledges support from NSF grant DMS-0071037}
 and Irena Swanson\unskip\footnote{**}{
\footfont The second author acknowledges support from NSF grant DMS-9970566}
}
\unskip\footnote{ }{{\footitfont 2000 Mathematics Subject Classification.}
13P10, 13A35
}
\unskip\footnote{ }{{\footitfont Key words and phrases.}
\footfont Frobenius powers, Gr\"obner bases, tight closure, binomial ideals.}

\vskip 4ex

\vskip 0.5cm

{\bgroup
\narrower 
\noindent {\bf Abstract.}
It is an open question whether
tight closure commutes with localization in quotients of a
polynomial ring in finitely many variables over a field.
Katzman~[K] showed that tight closure 
of ideals in these rings
commutes with
localization at one element
if for all ideals $I$ and $J$ in a polynomial ring
there is a linear upper bound in $q$ on
the degree in the least variable of reduced Gr\"obner bases
in reverse lexicographic ordering of the ideals of the form $J + \Fq I$.
Katzman conjectured that this property would always be satisfied.
In this paper we prove several cases of Katzman's conjecture.
We also provide 
an experimental analysis (with proofs) of 
asymptotic properties of Gr\"obner bases
connected with Katzman's conjectures.

\egroup
}
\bigskip

\section{Introduction}

Throughout this paper $F$ is a field of positive prime characteristic $p$,
$R$ is a finitely generated polynomial ring $F[x_1, \ldots, x_n]$ over $F$,
$J$ and $I$ denote ideals of $R$,
and $q = p^e$ denotes a power of $p$, where $e$ is a non-negative integer.
Then $\Fq I$ is the {\it $e\,$th Frobenius power of $I$},
defined by
$$
\Fq I := (i^q | i \in I).
$$
It follows that if $I$ is generated by $f_1, \ldots, f_r$,
then $\Fq I$ is generated by $f_1^q, \ldots, f_r^q$.

The main motivation for our work in this paper
is the theory of tight closure, in which
Frobenius powers of ideals play a central role.
In particular,
we address the question of whether tight closure commutes with localization.
The basics of tight closure can be found in the first few sections of~[HH];
however, in the following paper no knowledge of tight closure will be needed.

The polynomial ring $R$ is a regular ring,
so every ideal in $R$, and in the localizations of $R$,
is tightly closed [HH, Theorem 4.4],
and hence tight closure commutes with localization in $R$.
However, it is not known if
tight closure commutes with localization in quotient rings $R/J$ of $R$,
even for the special case of localization
at a multiplicatively closed set 
$\{1, r, r^2, r^3, \ldots\}$,
generated by one element $r \in R/J$.
Katzman~[K] showed that for this special case
it suffices to consider the case $r = x_n$ 
(by possibly modifying $R$, $I$, and $J$).
Katzman also proved that a positive answer to the question of tight closure
commuting with localization at $x_n$ would be provided by a positive answer
to the following conjecture.

\conj (Katzman~[K, Conjecture 4])
 Let $R=F[x_1,...,x_n]$
where $F$ is a field of characteristic $p$, and let
$I$ and $J$ be ideals of $R$.
Let $G_q$ be the reduced Gr\"obner basis for the
ideal $J + \Fq I$ with respect to the reverse lexicographic ordering.
Then there exists an integer $\alpha$
such that
the degrees in $x_n$ of the elements of $G_q$ are bounded above
by $\alpha q$.
\endb

The {\it (graded) reverse lexicographic ordering} on monomials in 
$x_1, \ldots, x_n$
is defined by
$x_1^{a_1}x_2^{a_2} \cdots x_n^{a_n} <x_1^{b_1}x_2^{b_2} \cdots x_n^{b_n}$
if $\sum_i a_i < \sum_i b_i$,
or if $\sum_i a_i = \sum_i b_i$
and $a_i > b_i$ for the last index $i$ at which $a_i$ and $b_i$ differ.
For background on reduced Gr\"obner bases, and Buchberger's
algorithm for finding these bases, see for example [CLO].

Katzman's conjecture holds trivially when $J = (0)$,
since Frobenius powers commute with sums in rings of characteristic $p$,
and hence the reduced
reverse lexicographic Gr\"obner basis for $\Fq I$
consists of the $q$th powers of elements of
the reduced Gr\"obner basis for $I$.
The other known cases are due to Katzman, 
who proved that the conjecture also holds whenever $J$ is generated
by monomials [K, Theorem 8],
and 
whenever $J$ is generated by binomials and simultaneously $I$ is
generated by monomials [K, Corollary 11].
There are classes of examples for which it is known
that tight closure commutes with localization
but for which Katzman's conjecture has not been proved; in particular,
one such class,
due to Smith [S],
consists of ideals $I$ and $J$ for which $J$ is a 
binomial ideal and $I$ is arbitrary.
Since the question of whether
tight closure commutes with localization
has so far defied proof for quotient rings of polynomial rings,
accordingly the proof of Katzman's conjecture is expected to be hard.
Difficulties in finding a general proof
include the dependence of Gr\"obner bases
on the characteristic of the field $F$
and the dependence of Gr\"obner bases
on raising a subset of the generators to powers.

In this paper we study the asymptotic behavior of
three functions of $q$ associated to the family
of reduced reverse lexicographic Gr\"obner bases $G_q$ for the ideals
$J + \Fq I$, namely

\item{(1)} the maximum of the $x_n$-degrees of the elements of $G_q$
(as in Katzman's conjecture),
also written as the $x_n$-degree of $G_q$ and denoted $\xq$;
\item{(2)} the maximum of the total degrees of the elements of $G_q$,
also referred to as the total degree of $G_q$ and denoted $\tq$; and
\item{(3)} the cardinality $\cq$ of $G_q$.
\par
\noindent Since for any ideals $I$ and $J$, $\xq \le \tq$
for all $q$, a linear upper bound for $\tq$ also implies
Katzman's conjecture.  In the special cases mentioned above for which
Katzman's conjecture is known to be true,
namely in which $J=(0)$, 
or $J$ is generated by monomials (with arbitrary ideal $I$), 
or $J$ is generated by binomials and $I$ by monomials,
Katzman's proof also shows a linear
upper bound for the function $\tq$ [K].

In Section~2 of this paper we prove (Theorem~2.1) %
that Katzman's conjecture
holds for polynomial rings in one or two variables with arbitrary
ideals $I$ and $J$, and find a linear upper bound for $\tq$ and
a constant upper bound for $\cq$ as well.
(As part of the proof of this theorem, we include a 
review of the steps of the Buchberger algorithm
for reduced Gr\"obner basis computation.)
In Sections 3 and 4 we provide further information about
the specific form of the functions $\xq$ and $\tq$, as
well as the function $\cq$, in the more
restrictive case in which $I$ and $J$ are both
principal binomial (and not monomial) ideals,
both to
gain better understanding of these functions and
to find (constructive) proofs of special cases of Katzman's conjecture
with potential for application in more general cases. 
In Section 3 we compute (in Theorem 3.2) Gr\"obner
bases for the ideals $J+\Fq I$ for ideals
$I=(x^u(x^v-gx^w))$ and $J=(x^a(x^b-hx^c))$ whenever
$g$ and $h$ are units, $\gcd(x^u,x^w)=1=\gcd(x^b,x^c)$, 
and $(x^v-gx^w,x^b-hx^c)=R$, and 
hence obtain a constructive proof of upper bounds
for $\xq$, $\tq$, and $\cq$ in this case.

In Theorem~3.3 we prove 
that for ``most'' principal
binomial ideals $I$ and $J$,
there is a change of variables that converts
$I$ and $J$ into monoidal ideals,
i.e.,
ideals generated by binomials of the form
$u-v$, where $u$ and $v$ are (monic) terms,
so that the coefficients are restricted to $+1$ and $-1$.
This change of variables preserves both 
the reverse lexicographic ordering on
the monomials and all three of the functions $\xq$, $\tq$, and $\cq$.
When $I$ and $J$ are monoidal ideals, the quotient rings
$R/(J+\Fq I)$ are monoid rings over $F$ for finitely presented
commutative monoids, and the Gr\"obner bases
for the ideals $J+\Fq I$
can also be considered to be finite complete rewriting systems
in the category of commutative monoids.

In Section~4 we study the asymptotic behavior of the three functions
$\xq$, $\tq$, and $\cq$ for constructions of
the reduced Gr\"obner bases $G_q$ for
a wide range of examples of principal monoidal ideals $I$ and $J$.
We give examples illustrating that
the three functions 
can be linear, periodic,
or have linear expressions holding only for $q$ sufficiently large;
in addition, we show examples in which
the cardinality and the $x_n$-degree
of the Gr\"obner bases can be bounded above by a constant for all $q$.
We also discuss 
the dependence of the three functions on the characteristic
$p$ of the field $F$ for several of the examples.
Section~4 ends with a
table summarizing the range of types of behavior of the Gr\"obner bases
we computed.
Finally, in the Appendix we include a sample of the
Macaulay2 [GS] code 
we used to generate Gr\"obner bases for small values of $q$
as an aid to our proofs.

\section{Katzman's conjecture for one and two variables}

In this section we prove the special case of Katzman's conjecture
for $n \le 2$.

\thm
\label{\thmdimonetwo}
Katzman's conjecture holds
when $R$ is a polynomial ring in one or two variables over $F$.
Moreover, 
for any ideals $I$ and $J$ in $R$ and
reduced Gr\"obner basis $G_q$ for the ideal 
$J + \Fq I$
with respect to the reverse lexicographic ordering,
there exist integers $\alpha$ and $\beta$ such that 
the $x_n$-degree and total degree functions satisfy
$\xq \le \tq \le \alpha q$ and the cardinality function satisfies
$\cq \le \beta$ for all $q$. 
\endb

\proof
If $R$ is a polynomial ring in one variable, then
$R$ is a principal ideal domain,
so $I = (f)$ and $J = (g)$ for some $f, g \in R$.
In this case $J + \Fq I$ is also a principal ideal,
and the reduced Gr\"obner basis of $J + \Fq I$
consists of the element $\gcd(g,f^q)$,
whose total degree is bounded above by $q \dg f$.
Then if we define $\alpha:= \dg f$ and $\beta:=1$, we obtain
$\xq \le \tq \le \alpha q$ and $\cq = \beta$ for all $q$.

Next suppose that $R$ is a polynomial ring in two variables
$x$ and $y$ over $F$.
By earlier observations, 
we may assume that $I$ and $J$ are non-zero ideals.
Let $S$ be a generating set for the ideal $J$,
and $T$ a generating set for $I$;
choose $S$ and $T$ so that the leading coefficients 
of all of their elements are $1$.
Define $T_q := \{t^q ~|~ t \in T\}$ to be the
corresponding generating set for $\Fq I$.

We apply the Buchberger algorithm with the reverse lexicographic ordering to
compute a Gr\"obner basis of $J + \Fq I$, starting with the generating set
$S \cup T_q$.
At each step, a partial Gr\"obner basis 
$B_{i-1}:=S \cup T_q \cup \{p_1,...,p_{i-1}\}$
has been found,
and an $S$-polynomial of a pair of elements in $B_{i-1}$ is
computed and reduced with respect to all of the elements in this basis.
If the result is non-zero,
the polynomial is divided by its leading coefficient
and the resulting monic polynomial
is added as the element $p_i$ to form the basis $B_{i}$.
When there are no non-zero reduced $S$-polynomials remaining,
this creates a Gr\"obner basis $B:=S \cup T_q \cup \{p_1,...,p_k\}$
for $J + \Fq I$ with respect to the reverse lexicographic ordering,
where each element of $B$ is a monic polynomial,
and for each $1 \le i \le k$, all of the terms of
the polynomial $p_i$ are reduced with respect to
$S \cup T_q \cup \{p_1,...,p_{i-1}\}$.

In order to compute the
reduced Gr\"obner basis $G_q$ of $J + \Fq I$,
we need to reduce the Gr\"obner basis $B$.
For each polynomial $r \in B$,
replace $r$ in the basis with the monic polynomial obtained by
reducing all of the terms of $r$ with respect
to the elements of $B \setminus \{r\}$, and dividing by the
resulting leading coefficient.
Repeat this process for all of the polynomials in the basis,
removing any zero polynomials that result, until no more
reduction can be done.
This gives the reduced Gr\"obner basis $G_q$ for 
$J + \Fq I$ [CLO, Prop. 2.7.6].

The total degree 
of the
reduced Gr\"obner basis $G_q$ for $J + \Fq I$
will be at most the total degree for the basis $B$.
To compute bounds on these degrees,
we first need to describe the polynomials $p_i$ more carefully.

Let $x^a y^b$ be the leading term of a non-zero element $p$ of 
$S$. 
In particular,  
since $J=(S) \ne (0)$, there is a non-zero monic polynomial
$p' \in J$, and
by adding the element $xyp' \in J$ to the set $S$ if necessary,
we may assume (for ease of notation) that both $a$ and $b$ are non-zero.
For each $1 \le i \le k$, let $x^{a_i}y^{b_i}$ be the
leading term of the polynomial $p_i$ in $B$.  Since $p_i$ is
reduced with respect to $S$, either $0 \le a_i<a$ or $0 \le b_i<b$,
or both.
If $i>j$, then $p_i$ is also reduced with respect to $p_j$.
More specifically, at each step of the algorithm described above,
when $p_i$ is computed, (at least)
one of four possible cases occurs.
Either
\item{(1)} $0 \le a_i<a$ and $a_i \ne a_j $ for all $1 \le j \le i-1$,
\item{(2)} $0 \le a_i<a$ and for some $j<i$, $a_i = a_j$ and $b_i<b_j$,
\item{(3)} $0 \le b_i<b$ and $b_i \ne b_j $ for all $1 \le j \le i-1$,
or
\item{(4)} $0 \le b_i<b$ and for some $j<i$, $b_i = b_j$ and $a_i<a_j$.

\noindent

In cases
(2) and (4), the total degree of $p_i$ is strictly less than
the maximal total degree of the previous
basis $S \cup T_q \cup \{p_1,...,p_{i-1}\}$.
In cases (1) and (3), the total degree of the polynomial
$p_i$, which is a reduction of an $S$-polynomial
of a pair of elements in the previous basis,
can be at most twice the maximal total degree of the previous basis
(by definition of S-polynomials).
Note that cases (1) and (3) can occur at most
$a+b$ times during the algorithm.
The maximal total degree of elements in $S \cup T_q$ satisfies
$$
\eqalign{
\dg(S \cup T_q)&=max\{\dg(S),\dg(T_q)\}
=max\{\dg(S),q \cdot \dg(T)\} \cr
&\le
q \cdot max\{\dg(S),\dg(T)\}.\cr}
$$
Thus the total degree of the basis $B$ is
at most $2^{a+b} \cdot q \cdot max\{\dg(S),\dg(T)\}$.
If we define the constant
$\alpha:=2^{a+b} \cdot max\{\dg(S),\dg(T)\}$, then
this proves that $\tq \le \alpha q$.
Since for all $q$, $\xq \le \tq \le \alpha q$,
therefore Katzman's conjecture holds in the case
in which the polynomial ring has two variables.

Finally, to get the bound on the cardinality of the
reduced Gr\"obner basis $G_q$ of $J+\Fq I$, 
note that although the element $p \in S$ with
leading term $x^ay^b$ may have been reduced or
removed in the reduction process to construct $G_q$ from $B$,
no polynomial that remains in $G_q$ may have
leading term divisible by $x^ay^b$.
For each number $0 \le a' < a$ and
$0 \le b' < b$, there
can be at most one polynomial in $G_q$ with leading term
of the form $x^{a'}y^{*}$ for any number $*$,
and at most one polynomial in $G_q$ with
leading term
$x^{*}y^{b'}$.  Then the cardinality of $G_q$ satisfies
$|G_q| \le a+b$.  By defining
the constant $\beta := a+b$, we obtain
$\cq \le \beta$.
\qed

\section{Principal binomial ideals:  General constructions}

For the remainder of the paper we direct our attention to the case in which
the ideals $I$ and $J$ are principal and binomial, and
obtain more detailed information about the specific form 
of the degree functions
$\xq$ and $\tq$, as well as the cardinality function $\cq$.
We begin by considering arbitrary monoidal
binomials which generate the whole ring.

\lemma
\label{\lmwholering}
Let $F$ be a field and
let $R = F[x_1, \ldots, x_n]$ be a polynomial ring in $n$ variables over $F$.
Let $x^v-\g x^w, x^b-\h x^c \in R$,
where $v,w,b,c$ are $n$-tuples of non-negative integers,
$\g, \h$ are units in $F$,
$\gcd(x^v, x^w) = 1 = \gcd(x^b, x^c)$,
and in reverse lexicographic ordering,
$x^v > x^w$ and $x^b > x^c$.
Assume that $(x^v-\g x^w, x^b-\h x^c) = R$.
Then $w = c = \underline{0}$,
and there is a positive rational
number $l$ such that $v_i=lb_i$ for all $i$.
\endb

\proof
If the conclusion holds after tensoring 
with the algebraic closure $\overline F$
of $F$ over $F$,
then it also holds in $R$.
So without loss of generality we may assume that $F$ is algebraically closed.

The hypothesis on the ordering implies that $v$ and $b$
are both non-zero.
If both $w$ and $c$ are also both non-zero,
then $R = (x^v-\g x^w, x^b-\h x^c) \subseteq (x_1, \ldots, x_n)R$,
which is a contradiction.
So either $w$ or $c$ is zero;
without loss of generality suppose that $w = \underline{0}$.

Choose any $(\k_1, \ldots, \k_n) \in F^n$ 
such that $\k^v = \g $.
Choose $i$ such that $v_i > 0$.
Since $v_i>0$, then $\k_i \not = 0$,
and $\k_i$ depends on the choices of the other $\k_j$ by the relation
$$
\k_i = \g ^{1/v_i} \prod_{j \not = i, v_j \not = 0} \k_j^{-v_j/v_i}
$$
(for some choice of the $v_i$th roots).
Since $\k^v-\g \k^w=\g-\g=0$, 
the assumption that $(x^v-\g x^w, x^b-\h x^c) = R$
implies that $\k^b- \h \k^c$ is a unit in $F$.
In particular,
for all indices $j$ with $v_j = 0$, 
any choice of $x_j=k_j \in F$ for these indices makes
$$
\k^b- \h \k^c=
\g^{b_i/v_i} \prod_{j \not = i, v_j \not = 0} \k_j^{b_j -b_i(v_j/v_i)}
\prod_{v_j = 0} x_j^{b_j}
- \h \g ^{c_i/v_i}
\prod_{j \not = i, v_j \not = 0} \k_j^{c_j -c_i(v_j/v_i)}
\prod_{v_j = 0} x_j^{c_j}
$$
a unit in $F$.

Suppose that $m$ is an index such that $m \not = i$ and $v_m = 0$.
If $b_m>0$ and $c_m>0$, then for the choice of $x_m=0$ 
the displayed expression above is $\k^b- \h \k^c=0$,
which is not a unit, giving a contradiction.
If $b_m=0$ and $c_m \not = 0$,
then with the choice of $k_j=1$ and $x_j=1$ for all $j \not = i,m$,
$k_i=g^{1/v_i}$, 
and $x_m=(h^{-1}g^{(b_i-c_i)/v_i})^{1/c_m}$,
the expression is again zero and not a unit, giving a contradiction.
Similar choices show that the case in which 
$b_m \not =0$ and $c_m = 0$ cannot occur. 
Therefore when $v_m = 0$, we have that $b_m=c_m=0$.
Thus $b_m -c_m = 0 = (b_i-c_i)(v_m/v_i)$
for all indices $m \not = i$ with $v_m = 0$.

Next let $m$ be any index such that $m \not = i$ and $v_m \not = 0$.
If in addition $\k_1, \ldots, \k_n$ are all chosen to be non-zero, then
$$
\k^{b-c}- \h=
\g^{(b_i-c_i)/v_i} \prod_{j \not = i} \k_j^{b_j -c_j -
(b_i-c_i)(v_j/v_i)}-\h
$$
is also a unit in $F$. 
If $b_m -c_m -
(b_i-c_i)(v_m/v_i) \not = 0$, then
for the choice
of $k_j=1$ for all $j \not = i,m$, and
the choice of 
$k_m=[hg^{-(b_i-c_i)/v_i}]^{1/(b_m -c_m -
(b_i-c_i)(v_m/v_i))}$, 
the expression above is $k^{b-c}-h=0$, giving a contradiction.
So $b_m -c_m - (b_i-c_i)(v_m/v_i) = 0$, or
$b_m -c_m = (b_i-c_i)(v_m/v_i)$, when $v_m \not = 0$ also.

Thus for all $j \not = i$, we have that
$b_j -c_j = (b_i-c_i)(v_j/v_i)$ and $v_j/v_i$ is non-negative.
By hypothesis $x^b > x^c$ in the reverse lexicographic ordering, so
we must have $b_i-c_i>0$ and 
$b_j \ge c_j$ for all $j$.
By the assumption that $\gcd(x^b,x^c) = 1$,
it follows that $c = \underline{0}$.
Then $b_iv_j=b_jv_i$
for all $j$, and since $v_i \ne 0$ and $b \ne \underline{0}$, 
$b_i \ne 0$ as well.
Therefore if we define 
the positive rational number $l:=v_i/b_i$, then $v_j=l b_j$ for all $j$.
\qed

This result leads to the following definition.
Two binomials $x^u(x^v-gx^w)$ and $x^a(x^b-hx^c)$ with
$x^v>x^w$ and $x^b>x^c$ are
of the {\it same type} if
there are non-negative integers $l$ and $m$ and
$n$-tuples $B$ and $C$ of non-negative integers
with
$x^B > x^C$
such that $v = lB$, $w = lC$, $b = mB$, and $c = mC$;
in this case, we say the binomials are of {\it type} $(B,C)$. 
With this notation the Lemma above says
that if the ideal generated by two non-monomial binomials
is the whole ring, then the two binomials
are both of type $(B, (0,\ldots,0))$ for some $B$,
and neither binomial is a multiple of any variable.
The corresponding result fails for a $3$-generated binomial ideal;
for example,
the three binomials $x_1-1$, $x_2-1$, $x_1 x_2 -2$
generate the whole ring,
yet no two of the three binomials are of the same type.

The following theorem shows that
for the ideals considered in Lemma~\lmwholering,
one can bound the number of elements in the reduced Gr\"obner bases,
as well as give constructive upper bounds for the $x_n$-degree and total degree.

\thm
\label{\prrelprim}
Let $F$ be a field of positive prime characteristic $p$
and $R = F[x_1, \ldots, x_n]$ a polynomial ring in $n$ 
variables over $F$.
Let $I = (x^u(x^v-\g x^w))$
and $J = (x^a(x^b-\h x^c))$ be ideals in $R$,
where $u,v,w,a,b,c$ are $n$-tuples of non-negative integers,
$\g, \h$ are units in $F$,
$\gcd(x^v, x^w) = 1 = \gcd(x^b, x^c)$,
and in reverse lexicographic ordering,
$x^v > x^w$ and $x^b > x^c$.
Assume that $(x^v-\g x^w, x^b-\h x^c) = R$.
Then for $q$ sufficiently large, the maximal $x_n$-degree
of the Gr\"obner basis of $J+ \Fq I$ is
$\delta(q)\le \max((u_n+v_n)q, a_n+b_n)$, the maximal total degree
is
$\Delta(q)\le \max((|u|+|v|)q, |a|+|b|)$,
and
the cardinality of the Gr\"obner basis is $\cq \le 4$.
\endb

\proof
By Lemma~\lmwholering, $w=c=0$ and the 
generators of $I$ and $J$ have the same type.
Then $\Fq I=(x^{qu} (x^{qv}-g^q))$ and $J=(x^a(x^b-h))$.
We will explicitly compute a Gr\"obner basis for $J+ \Fq I$.

The hypothesis that $(x^v-\g, x^b-\h) = R$
implies that there are polynomials $r,s \in R$ with
$r(x^v-\g x^w)+s(x^b-\h x^c)=1$.  
Taking $q$th powers
of both sides and then multiplying by $\lcm(x^{qu},x^a)$
yields
$$r^q{\lcm(x^{qu},x^a) \over x^{qu}}x^{qu}(x^{qv}-\g^q)+
[s^q(x^b-\h x^c)^{q-1}]{\lcm(x^{qu},x^a) \over x^{a}}x^{a}(x^b-\h x^c)=
\lcm(x^{qu},x^a).
$$
Thus $J + \Fq I$ contains $\lcm(x^{qu},x^a)$.
Computation of the S-polynomials of this monomial
with the two generators 
of $J + \Fq I$ shows that
$$
{1 \over g^q}S(x^{q(u+v)}-g^qx^{qu},\lcm(x^{qu},x^a))=
{1 \over g^q}{\lcm(x^{q(u+v)},x^a) \over x^{q(u+v)}} g^qx^{qu} =
{\lcm(x^{q(u+v)},x^a) \over x^{qv}}
\hbox{\ \ and \ \ }
$$
$$
{1 \over h}S(x^{a+b}-hx^{a},\lcm(x^{qu},x^a))=
{1 \over h}{\lcm(x^{a+b},x^{qu}) \over x^{a+b}} hx^a =
{\lcm(x^{a+b},x^{qu}) \over x^b}
$$
are also in $J + \Fq I$.

Let $E_j := {\lcm(x^{a+jb},x^{qu}) \over x^{jb}}$.
By the $S$-polynomial calculation above, $E_1 \in J + \Fq I$.
If $E_j \in J + \Fq I$,
then so is
$$
{1 \over \h} S(x^a(x^b-\h), E_j)
=
{\lcm(x^{a+b}, {\lcm(x^{a+jb},x^{qu}) \over x^{jb}}) \over x^b}.
$$
The exponent of $x_i$ in this equals
$\max(a_i+b_i, \max(a_i+jb_i,qu_i)-jb_i)-b_i$
$=\max(a_i, \max(a_i-b_i,qu_i-(j+1)b_i))$
$=\max(a_i, qu_i-(j+1)b_i)$,
which is the same as the exponent of $x_i$ in $E_{j+1}$.
Thus the monic S-polynomial above is
${1 \over \h} S(x^a(x^b-\h), E_j)=E_{j+1}$.
Therefore all the $E_j$ are in $J + \Fq I$.
Note that the exponent $\max(a_i,qu_i-jb_i)$ of
$x_i$ in $E_j$ is at least as large as the
exponent of $x_i$ in $E_{j+1}$ for all $i$, so
$E_j$ is a multiple of $E_{j+1}$ for each $j$.
Thus for sufficiently large $j$,
$E_j = E_{j+1}=E_{j+2}= \cdots$, and we 
denote this eventual monomial as $E_{\infty}$.
All of the $E_j$ are multiples of $E_{\infty}$.

Define the set
$$
B:=\left\lbrace
x^{qu}(x^{qv}-\g ^q),
x^a(x^b - \h),
{\lcm(x^{q(u+v)},x^a) \over x^{qv}}, E_{\infty} \right\rbrace;
$$
then $B$ is a basis of $J + \Fq I$.
The S-polynomial of the first two elements is
${\lcm(x^{q(u+v)},x^{a+b}) \over x^{qv}} \g ^q
-{\lcm(x^{q(u+v)},x^{a+b}) \over x^b}\h$,
which reduces modulo the third element in $B$ and modulo $E_1$
(i.e. modulo $E_{\infty}$) to zero.
The S-polynomial of the first and the third elements in $B$ is
$$
S\left(x^{qu}(x^{qv}-\g ^q), {\lcm(x^{q(u+v)},x^a) \over x^{qv}}\right)
=
{\lcm(x^{q(u+v)},{\lcm(x^{q(u+v)},x^a) \over x^{qv}}) \over x^{qv}}
\g ^q.
$$
The exponent of $x_i$ in this equals
$$
\max(qu_i + qv_i, \max(qu_i + qv_i, a_i) - qv_i) - qv_i
= \max(qu_i, a_i- 2 qv_i).
$$
For sufficiently large $q$,
if $v_i \not = 0$ then $\max(qu_i, a_i- 2 qv_i)=qu_i=\max(qu_i, a_i- qv_i)$,
and if $v_i=0$ then
$ \max(qu_i, a_i- 2 qv_i)=
\max(qu_i, a_i)=
\max(qu_i, a_i- qv_i)$.
Since $\max(qu_i, a_i- qv_i)$ also equals the exponent of $x_i$
in the third element of the basis $B$, this shows that
the S-polynomial of the first and the third element of $B$
reduces to $0$.
The S-polynomial of the first and the fourth elements in $B$ is
$$
S\left(x^{qu}(x^{qv}-\g ^q), E_{\infty}) \right)
=
{\lcm(x^{q(u+v)},E_{\infty}) \over x^{qv}}\g ^q.
$$
The exponent of $x_i$ in this equals, for $j$ sufficiently large,
$$
\max(qu_i, \max(a_i - qv_i, qu_i - j b_i - qv_i))
= \max(qu_i, a_i- qv_i),
$$
which is the same as the exponent of $x_i$
in the third element of $B$.
Thus the S-polynomial of the first element of $B$ 
with any other element of $B$ reduces to 0.
The S-polynomial of the second and third elements is the monomial
$$
{\lcm\left(x^{a+b},{\lcm(x^{q(u+v)},x^a) \over x^{qv}}\right)
\over x^b} \h,
$$
for which the exponent of $x_i$ is
$\max(a_i, \max(qu_i - b_i, a_i-qv_i-b_i))
=\max(a_i, qu_i - b_i)$,
so that this S-polynomial is a multiple of $E_1$ and 
thus of $E_{\infty}$, and hence reduces to zero.
We have previously established that the S-polynomial
of the second and the fourth elements reduces to 0 modulo the given basis.
The last two elements of the basis $B$ are both monomials,
so their S-polynomial is 0 as well.
This proves that for $q$ sufficiently large the set $B$ is a
Gr\"obner basis of $J + \Fq I$ with respect to the
reverse lexicographic ordering.  

Although the Gr\"obner basis $B$ may not be reduced,
the reduced reverse lexicographic
Gr\"obner basis $G_q$ for $J+\Fq I$ will have cardinality
and degrees at most those of $B$.
Thus we can read off upper bounds for
the three functions for $q$ sufficiently large, and
find that
$\delta(q)\le \max((u_n+v_n)q, a_n+b_n)$,
$\Delta(q)\le \max((|u|+|v|)q, |a|+|b|)$, and
$c(q) \le 4$.
\qed

Next we use Lemma~\lmwholering~to show
that the principal monoidal
ideals cover ``most'' of the possibilities for
principal binomial ideals.

\thm
\label{\indcoeftwo}
For any principal binomial (non-monomial) ideals $I$ and $J$
which are generated by binomials that are not of the same type,
there is a change of variables 
under which $I$ and $J$ become principal monoidal ideals.
Furthermore,
this change of variables preserves the 
reverse lexicographic ordering and the
three functions 
$\delta(q)$, $\Delta(q)$, and $c(q)$.
\endb

\proof
Let $F$ be a field of positive prime characteristic $p$
and $R = F[x_1, \ldots, x_n]$ a polynomial ring in $n$ variables over $F$.
Since Gr\"obner bases are unchanged if we pass to 
$\overline F[x_1, \ldots, x_n]$,
where $\overline F$ is the algebraic closure of $F$,
without loss of generality we may assume that $F$ is 
algebraically closed.

Let $I$ and $J$ be arbitrary principal binomial (non-monomial) ideals
that are not of the same type.
We can write $I = (x^u(x^v-\g x^w))$
and $J = (x^a(x^b-\h x^c))$,
where $u,v,w,a,b$, and $c$ are $n$-tuples of non-negative integers,
$\g$ and $\h$ are units in $F$, $x^v>x^w$ and $x^b>x^c$ in the
reverse lexicographic ordering,
and $\gcd(x^v, x^w) = 1 = \gcd(x^b, x^c)$.

Case I.
Suppose 
there exist non-zero elements $\k_1$,
$\ldots$,
$\k_n$ in $F$
such that $\k^v-\g \k^w = 0 = \k^b-\h \k^c$. 
In this case
under the variable change $x_i \mapsto \k_i x_i$ for all $i$,
the reverse lexicographic ordering is preserved,
and the generator of the image $\widetilde I$
of $I$ under this ring automorphism is
$\k^ux^u (\k^vx^v-\g \k^w x^w)$.
After dividing through by the non-zero element
$\k^u \k^v= \k^u \g \k^w$ of $F$,
this generator becomes $x^u (x^v-x^w)$.
A similar computation holds for the generator of
the image $\widetilde J$ of $J$; hence
the generators of $\widetilde I$ and $\widetilde J$ are monoidal.
As this ring automorphism preserves the reverse lexicographic ordering,
it maps Gr\"obner bases to Gr\"obner bases.
Since this
change of variables is linear, the functions $\xq$, $\tq$, and $\cq$
will also be preserved.

Case II.
Suppose that 
there do not exist non-zero elements $\k_1$, $\ldots$, $\k_n$ in $F$
such that $\k^v-\g \k^w = 0 = \k^b-\h \k^c$.

Case IIa.
Suppose Case II holds and also that 
$v_i + w_i > 0$ and $b_i + c_i = 0$ for some index $i$.
There is another index $j$ for which either $b_j > 0$ or $c_j > 0$,
but not both, since $x^b>x^c$ and $\gcd(x^b, x^c)=1$.
By performing the change of variables $x_j \mapsto \h^{1/b_j}$
(respectively $x_j \mapsto (\h^{-1})^{1/c_j}$) and $x_m \mapsto x_m$
for all $m \ne j$,
the generator $x^a(x^b-\h x^c)$ of $J$ is mapped to
a scalar multiple of $x^a(x^b-\tilde \h x^c)=x^a(x^b-x^c)$
with unit $\tilde \h=1$.
At the same time, the generator of $I$
changes to a scalar multiple of $x^u(x^v-\tilde \g x^w)$
for 
another unit $\tilde \g$ in $F$.
Since either $v_i > 0$ or $w_i > 0$, we can similarly
replace $x_i$ by an appropriate scalar multiple of itself
so that $x^u(x^v-\tilde \g x^w)$ is mapped to a scalar multiple of
$x^u(x^v- x^w)$.
Since $b_i = c_i = 0$,
the unit $\tilde \h = 1$ remains unchanged under this second map.
As in Case I, this change of variables preserves the ordering 
and the three functions associated to the Gr\"obner bases.

Case IIb.
Suppose Case II holds and 
$v_i + w_i = 0$ and $b_i + c_i > 0$ for some index $i$.
An argument similar to Case IIa shows this case as well.

Case IIc.
Suppose Case II holds and that for all indices $i$,
$v_i + w_i > 0$ if and only if $b_i + c_i > 0$.
Let $T$ be the set of indices $m$ for which $v_m>0$,
let $U$ be the set of indices $m$ for which $w_m>0$,
and let $S:=T \cup U$.
Let 
$$b_+ := \cases{b_j&if $j \in T$\cr 0&if $j \not \in T$\cr}
\hskip 0.5in and \hskip 0.5in
b_- := \cases{b_j&if $j \in U$\cr 0&if $j \not \in U$,\cr}
$$
and define $c_+$ and $c_-$ similarly.  Then 
$b = b_+ + b_-$ and $c = c_+ + c_-$.

Define new variables $y_m$ over $F$,
where $m$ varies over the set $S$.
We will denote the restrictions of the tuples
$v$, $w$, $b_+$, $b_-$, $c_+$, and $c_-$ to tuples in the indices
of $S$ by the same notation.
Consider the ideal
$(y^{v+w}-\g, y^{b_+ + c_-} -\h y^{b_- + c_+} )$
in $F[y_m | m \in S]$.

Since the non-leading (monic) term of the
first generator is $1$, it follows directly that
$y^{v+w}>1$ and 
$\gcd(y^{v+w},1)=1$.
For the second generator, the indices $m$ for which
$(b_+)_m>0$ satisfy that both $(b_-)_m =0= (c_+)_m$, and
similarly for $c_-$, so the supports of the two terms are
disjoint.  Then $\gcd(y^{b_+ + c_-},y^{b_- + c_+})=1$ and
either $y^{b_+ + c_-}>y^{b_- + c_+}$ or 
$y^{b_- + c_+}>y^{b_+ + c_-}$.

Suppose that $\underline{\tilde k}$ is a tuple with entries in $F$
(and indices in $S$) for which
$\tilde k^{v+w}-\g=0= \tilde k^{b_+ + c_-} -\h \tilde k^{b_- + c_+}$.
Since the product of all of the $\tilde k_m$ divides $\tilde k^{v+w}$, 
the first equation shows that all of the entries of 
$\underline{\tilde k}$ are non-zero.  
Define the $n$-tuple $\underline{k} \in F^n$ by
$k_j:=\tilde k_j$ for $j \in T$, $k_j:=\tilde k_j^{-1}$ for $j \in U$,
and $k_j:=1$ for $j \not \in S$.  Then $k_1, \ldots, k_n$
are non-zero elements in $F$ for which
$0=k^w(\tilde k^{v+w}-\g)=k^w(k^{v-w}-\g)=k^v-\g k^w$
and $0=k^{b_-+c_-}(\tilde k^{b_+ + c_-} -\h \tilde k^{b_- + c_+})=
k^{b_-+c_-}(k^{b_+ - c_-} -\h k^{-b_- + c_+})=k^b-\h k^c$,
contradicting the hypothesis of Case II.
Therefore the equations 
$y^{v+w}-\g=0=y^{b_+ + c_-} -\h y^{b_- + c_+}$ have no solutions
over $F$.
Then Hilbert's Nullstellensatz says that
$(y^{v+w}-\g, y^{b_+ + c_-} -\h y^{b_- + c_+} )
=F[y_m | m \in S]$.

Applying Lemma~\lmwholering, we get that
either $b_+ + c_- = 0$ or $b_-+c_+ = 0$, and we
can write $y^{b_+ + c_-} -\h y^{b_- + c_+}$ as 
a scalar multiple of $y^{\hat b+\hat c}-\hat \h$
where $\hat b+\hat c$ is either $b_+ + c_-$ or $b_-+c_+$, and
$\hat h$ is $\h$ or $\h^{-1}$, respectively.
The last conclusion of Lemma~\lmwholering\ says
there is a positive rational number $l$ such that
$v+w = {l} (\hat b + \hat c)$.
If $b_+ + c_- = 0$, then $b = b_-$, $c = c_+$, and
$v+w=l(b_-+c_+)$, so $v=lc_+$ and $w=lb_-$, which
contradicts the assumption that both $x^v>x^w$ and $x^b>x^c$.
Therefore $b_-+c_+ = 0$, so $b = b_+$, $c = c_-$,
$v=lb$, and $w=lc$.
Therefore the generator $x^u(x^v-\g x^w)$ of the ideal $I$
is of the same type as the generator $x^a(x^b-\h x^c)$ of $J$.
But this contradicts the hypothesis that the generators
of $I$ and $J$ are of distinct types, so Case IIc cannot occur.
\qed

Motivated by the preceding theorem, for the remainder of
the paper we consider the case in which the ideals $I$ and $J$ are
principal and monoidal.

\section{Principal monoidal ideals:  Examples}

In this section we report on our calculations of reduced 
reverse lexicographic Gr\"obner bases,
together with the functions $\xq$,
$\tq$, and $\cq$,
for ideals of the form $J + \Fq I$,
where $I$ and $J$ are fixed principal monoidal ideals
and $q$ varies over powers of the
characteristic of the base field $F$.
In every example, the three functions
either are eventually (for $q>>0$) linear or constant functions,
or else eventually vary periodically between linear 
or constant functions.
For several of the examples, we also explore in more detail
the dependence of the three functions on the characteristic
$p$ of the field $F$.
The examples included in this Section were chosen
from among all of our computations
to illustrate all of the possible behaviors 
we observed for the three functions.

In the process of finding each of the following examples,
we used the symbolic computer algebra program Macaulay2 [GS]
to generate Gr\"obner bases for
ideals $J + \Fq I$ for small values of $q$ (usually three or four values),
and studied the patterns in these bases
to guide us in proving the structure of the Gr\"obner bases
for all values of $q$.
A sample of the Macaulay2 code used in our calculations is
provided in the Appendix.

We begin with an example in which the degree functions are
linear functions and the cardinality is a constant.

\prop
\label{\lmbdunb}
Let $R = \bbZ/3\bbZ[x,y,z]$,
$I = (y^2z-x^2)$, $J = (y^3-xy)$, $p=3$, and $q=3^e$.
Then the Gr\"obner basis of $J + \Fq I$
with respect to the reverse lexicographic ordering
(with $x_1=x$, $x_2=y$, and $x_3=z$, so that $z<y<x$)
is
$$
\lbrace
y^3-xy, x^{q-1}y^2z^q-x^{2q},
x^{2q}y-x^qyz^q,
x^{3q+1}-x^{2q+1}z^q
\rbrace.
$$
Therefore the maximal $z$-degree of the Gr\"obner basis elements
for $J + I^{[q]}$ is $\xq=q$,
the maximal total degree of the elements is $\tq=3q + 1$,
and the number of elements in the Gr\"obner basis is $\cq=4$
for all $q$.
\endb

\proof
Define $f:=y^3-xy$ and $g:=y^{2q}z^q-x^{2q}$,
so that $f$ and $g$ generate $J$ and $\Fq I$,
respectively.
Before computing $S$-polynomials,
we reduce $g$ modulo $(y^3-xy)$.
Note that for any monomial $x^ay^bz^c$ with
$b \ge 3$, the monomial reduces to
$x^{a+1}y^{b-2}z^c$.
Then the normal form of $x^ay^bz^c$ modulo $f$ is $x^{a+k}y^{b-2k}z^c$,
where $b-2(k-1) \ge 3$ and $b-2k<3$;
that is, $(b-3)/2 < k \le (b-1)/2$.
Then to find the normal form
for $y^{2q}z^q$, where $b=2q$, we
need $q-{3 \over 2} < k \le q-{1 \over 2}$, so
$k=q-1$, and the normal form is
$x^{q-1}y^2z^q$.
Therefore the polynomial $g$ reduces to $g':=x^{q-1}y^2z^q-x^{2q}$.

The polynomials $f$ and $g'$ are a basis for $J+\Fq I$.
Let $h$ denote their S-polynomial
$$
h:=S(f,g')=x^{q-1}z^qf-yg' = -x^qyz^q+x^{2q}y.
$$
The $S$-polynomial
$$
\eqalign{
S(g',h) &= x^{q+1}g'-yz^qh=-x^{3q+1}+x^qy^2z^{2q} \cr
 & \equiv -x^{3q+1}+x^{2q+1}z^q,\cr
}
$$
where $\equiv$ denotes a reduction using $g'$; let
$i:=x^{3q+1}-x^{2q+1}z^q$ denote the monic scalar multiple
of this polynomial.  
All of the remaining $S$-polynomials 
in the basis $\{f,g',h,i\}$
reduce to $0$.
Therefore the four elements indeed generate a Gr\"obner basis,
and since no element of the basis may be reduced by any other,
this Gr\"obner basis is also reduced.
This proves
that the 
maximal $z$-degree is of the elements of the Gr\"obner basis $\xq=q$, 
the maximal total degree is
$\tq=3q +1$, and
the cardinality is $\cq = 4$.
\qed

\note
Let $R = \bbZ/p\bbZ[x,y,z]$,
with $x,y,z$ variables over $\bbZ/p\bbZ$, where $p$ is
any prime and
$q$ varies over powers of $p$.
Let $I = (y^2z-x^2)$ and $J = (y^3-xy)$ be the same ideals as
in the example above.
In this case,
the same sets as in Proposition~\lmbdunb\ above are the 
reduced Gr\"obner bases
of the ideals $J + \Fq I$ in characteristic $p$ also.
Indeed, the proof above applies, since the hypothesis
that $p=3$ was never used in the proof.
\qed

The number of elements in the Gr\"obner bases
need not remain constant, as we prove next with
the ideals $I$ and $J$ from Proposition~\lmbdunb,
but with their roles switched.

\prop
\label{\lmbdunbtwo}
Let $R = \bbZ/3\bbZ[x,y,z]$, $I = (y^2z-x^2)$, $J = (y^3-xy)$,
and $q=3^e$.
Then the Gr\"obner basis of $I + \Fq J$
(roles of $I$ and $J$ exchanged)
with respect to the reverse lexicographic ordering
(with $z<y<x$)
is
$$
\eqalign{
\lbrace &
y^2z-x^2, y^{3q}-x^qy^q,
x^{2k}y^{3q-2k} - x^{q+2k}y^{q-2k} ~|~
1 \le k \le (q-1)/2\rbrace \cr
&\cup \lbrace x^{q-1+2j}y^{2q+1-2j} - x^{2q-1}yz^j ~|~ 1 \le j \le q\rbrace
\cup \lbrace x^{3q+1} - x^{2q+1}z^q \rbrace. \cr
}
$$
The corresponding functions for these ideals are
$\xq=q$, $\tq=3q + 1$
and $\cq = (3q+5) / 2$ for all $q$.
\endb

\proof
Define the polynomials
$f:=y^2z-x^2$, $g:=y^{3q}-x^qy^q$,
$h_k:=x^{2k}y^{3q-2k} - x^{q+2k}y^{q-2k}$ when
$1 \le k \le (q-1)/2$,
$r_j:=x^{q-1+2j}y^{2q+1-2j} - x^{2q-1}yz^j$ when
$1 \le j \le q$, and
$s:=x^{3q+1} - x^{2q+1}z^q$.
Since $q=3^e$, $q$ is odd, so $(q-1)/2$ is
an integer for all values of $e$.

Note that if $q=1$, there are no elements of the form $h_k$.
In this case, the Gr\"obner basis is
already included in the proof of Proposition \lmbdunb.

Next assume that $q>1$.
In this example each of the generators of both $I$ and $\Fq J$
is in normal form with respect to the other,
giving the first two elements $f$ and $g$ of the 
basis.
The $S$-polynomial
$$
S(f,g)=y^{3q-2}f-zg=-x^2y^{3q-2}+x^qy^qz
\equiv -x^2y^{3q-2}+x^{q+2}y^{q-2}=-h_1,
$$
where $\equiv$ denotes a reduction using $f$.
Repeating this for $1 \le k \le (q-3)/2$,
we get
$$
\eqalign{
S(f,h_k) &= x^{2k}y^{3q-2k-2} f - z h_k %
  = -x^{2k+2}y^{3q-2k-2} + x^{q+2k}y^{q-2k}z \cr
  &\equiv -x^{2(k+1)}y^{3q-2(k+1)} + x^{q+2(k+1)}y^{q-2(k+1)}
  = -h_{k+1}, \cr
}$$
where $\equiv$ denotes a reduction of the second
term using $f$.
Note that in this $S$-polynomial computation, we
required that the first $y$-exponent ${3q-2k-2} \ge 0$,
and to do the later reduction by $f$,
we needed that $y^2$ divides $y^{q-2k}$.
Then $3q-2k \ge 2$ and $q-2k \ge 2$,
so the first inequality is redundant, and the
second inequality says $k \le (q-2)/2$.
Since in this
Proposition we are assuming that $p=3$, so $q=3^e$
is always odd, the largest value that
$k$ can actually reach in this $S$-polynomial
computation is $(q-3)/2$.
Then the largest value
of $k$ for which a basis element $h_k$ is
produced is $(q-1)/2$.
Thus 
the entire 
set of elements $h_k$ is generated in the Buchberger algorithm.

The last element generated this way is
$h_{(q-1)/2}=x^{q-1}y^{2q+1} - x^{2q-1}y$.
Computing the $S$-polynomial of this and $f$ gives $S(f,h_{(q-1)/2})=-r_1$.
Again computing $S$-polynomials inductively
for $1 \le j \le q-1$, we get
$S(f,r_j)=-r_{j+1}$.
The last element generated in this latter step
is $r_q= x^{3q-1}y - x^{2q-1}yz^q$.

Finally, the $S$-polynomial
$S(f,r_q)$ reduces (using $f$) to the polynomial $-s$,
resulting in the last element
in the list of the 
basis elements.
It is straightforward to
check that with these basis elements all remaining $S$-polynomials
reduce to 0, hence the set is a Gr\"obner basis, 
and that the Gr\"obner basis is reduced.  The results
on the three functions then follow directly.
\qed

\note
Let $R = \bbZ/2\bbZ[x,y,z]$, so that the characteristic is $p=2$, and
let $I = (y^2z-x^2)$ and $J = (y^3-xy)$ be the same ideals as
in the example above.
In the proof above, in the computation of
the $S$-polynomials $S(f,h_k)$, we noted that
the number of polynomials of the form $h_k$ produced
satisfies $k \le (q-2)/2$.
When the characteristic
$p$ is even, then, the Gr\"obner basis computation can differ
from the proof above
at that point.
In fact, a proof very similar to the one above shows
that for $p=2$ the reduced Gr\"obner basis of $I + \Fq J$ is
$$
\{ y^2z - x^2, x^{2k}y^{3q-2k}-x^{q+2k}y^{q-2k},
x^{q+2j}y^{2q-2j}-x^{2q}z^j ~|~ 0 \le k \le (q-2)/2, 0 \le j \le q\}
$$
when $q > 1$.
Then the functions $\tq=3q$ and $\cq= {3 \over 2}q+2$
for $q > 1$ associated to these Gr\"obner bases
differ from the functions $\tq$ and $\cq$ computed
in Proposition \lmbdunbtwo\ with $p=3$.
Thus, not surprisingly,
the reduced Gr\"obner bases 
do depend on the characteristic of the
underlying field in general.
In this example, though, the $x_n$-degree $\xq=q$ is the same function
in both characteristics. 
\qed

In part (a) of the next Proposition, we show that
the function $\xq$ also can equal a constant.
In the Propositions above, the functions $\xq$ and $\tq$
are exactly equal to linear functions, and $\cq$ equals
either a linear or constant function, for all $q$.
As mentioned earlier, these
functions are not always this regular.
Part (b) of the next Proposition
illustrates functions $\tq$ and $\cq$
which are polynomials eventually but not at the start.

\prop
\label{\lmbdbd}
Let $R = \bbZ/2\bbZ[x,y,z]$,
$I = (x^2-y^2)$, $J = (xy-z^2)$ and $q=2^e$.
With the reverse lexicographic ordering (with $z<y<x$),

(a) the
reduced Gr\"obner basis for
$J + \Fq I$ is
$$
\{xy-z^2,x^{2q}-y^{2q},y^{2q+1}-x^{2q-1}z^2\},
$$
so that $\xq = 2$, $\tq=2q+1$ and $\cq = 3$ for all $q$, and

(b) the reduced Gr\"obner basis
for $I + \Fq J$ with $q \ge 2$ is
$$
\{x^2-y^2, y^{2q}-z^{2q}\},
$$
so that in this case
$\xq=2q$,
$$
{
\tq = \cases{3&if $q=1$, \cr
                           2q&if $q  \ge 2$ \cr}
}
\qquad {\rm and} \qquad
{
\cq = \cases{3&if $q=1$, \cr
                           2&if $q  \ge 2$.  \cr}
}
$$
for all $q$.
\endb

\proof
In the first part,
the S-polynomial of the generators $xy-z^2$ and $x^{2q}-y^{2q}$
of $J$ and $\Fq I$, respectively,
is $y^{2q+1}-x^{2q-1}z^2$.
The S-polynomials of these three polynomials all reduce to zero.
This verifies part (a). 
In part (b), when $q \ge 2$
the generator $x^qy^q-z^{2q}$ of $\Fq J$ reduces modulo $I$ to
$y^{2q}-z^{2q}$.
Since the leading terms $x^2$ and $y^{2q}$ of
$x^2-y^2$ and $y^{2q}-z^{2q}$ have no common factors,
their S-polynomial is zero.
\qed

The next example
shows that the function $\xq$ may also be
a function that is eventually
linear but not for small $q$.

\prop
\label{\lmbdperiodtwo}
Let $R = \bbZ/3 \bbZ [x,y,z,w]$,
$I = (x^5y^2zw-xy^3z^2w)$,
$J = (xy^2z^3w^2-x^3yzw^3)$  and $q=3^e$.
Then with the reverse lexicographic ordering (with $w < z < y < x$), the
reduced Gr\"obner basis of $J+I$ is
$$
\{xy^2z^3w^2-x^3yzw^3, x^5y^2zw-xy^3z^2w, x^7yzw^3-x^3y^2z^2w^3\},
$$
and for $q \ge 3$,
the
reduced Gr\"obner basis of $J + \Fq I$
is
$$
\lbrace
xy^2z^3w^2-x^3yzw^3,
x^{6q-1+2i}y^{{3q+1 \over 2}-i}zw^{{3q-1 \over 2}+i}
- x^{3q-2+2i}y^{2q+1-i}z^2w^{2q-1+i} |~
0 \le i \le \hbox{${3q-1 \over 2}$}\rbrace.
$$
Thus 
$$
{
\xq = \cases{3&if $q=1$, \cr
                           {7q-3 \over 2}&if $q  \ge 3$, \cr}
}
\qquad \qquad
{
\tq = \cases{12&if $q=1$, \cr
                           12q-1&if $q  \ge 3$  \cr}
}
$$
and 
$\cq={3q+3 \over 2}$ for all $q$.
\endb

\proof
The case $q = 1$ can be computed directly by Macaulay2 and is
left to the reader.

Now assume that $q \ge 3$.
Define $f := xy^2z^3w^2-x^3yzw^3$,
$i_0 := {3q-1 \over 2}$,
and 
$$
g_i := x^{6q-1+2i}y^{{3q+1 \over 2}-i}zw^{{3q-1 \over 2}+i}
- x^{3q-2+2i}y^{2q+1-i}z^2w^{2q-1+i}
$$
for $0 \le i \le i_0$.
In particular,
$g_{i_0} = x^{9q-2}yzw^{3q-2} - x^{6q-3}y^{q+3\over 2}z^2w^{7q-3 \over 2}$.

Observe that whenever $a \ge 1$, $b \ge 2$, $c \ge 3$ and $d \ge 2$,
the monomial $x^a y^b z^c w^d$ reduces to
$x^{a+2} y^{b-1} z^{c-2} w^{d+1}$ modulo $J$,
so that the normal form of the monomial
$x^a y^b z^c w^d$ modulo $J$
is $x^{a+2k} y^{b-k} z^{c-2k} w^{d+k}$,
where $k$ is the largest integer such that
$b-(k-1) \ge 2$ and $c-2(k-1) \ge 3$, i.e.,
$b \ge k+1$ and $c \ge 2k +1$.
In particular,
as $q \ge 3$,
the generator $x^{5q}y^{2q}z^qw^q-x^qy^{3q}z^{2q}w^q$ of $\Fq I$
reduces to
$g_0 = x^{6q-1}y^{3q+1 \over 2}zw^{3q-1 \over 2} - 
 x^{3q-2}y^{2q+1}z^2w^{2q-1}$.

Computing the $S$-polynomials, we begin with $S(f,g_i)$ for $i < i_0$.
$$
\eqalignno{
S(f,g_i) &=
x^{6q-2+2i} y^{{3q-3 \over 2}-i} w^{{3q-5 \over 2}+i} f - z^2 g_i, \cr
&= - x^{6q-1+2(i+1)} y^{{3q+1 \over 2}-(i+1)} z w^{{3q-1 \over 2}+(i+1)}
+ x^{3q-2+2i}y^{2q+1-i}z^4w^{2q-1+i} \cr
&\equiv - x^{6q-1+2(i+1)} y^{{3q+1 \over 2}-(i+1)} z w^{{3q-1 \over 2}+(i+1)}
+ x^{3q-2+2(i+1)}y^{2q+1-(i+1)}z^2w^{2q-1+(i+1)} \cr
&= -g_{i+1}, \cr
}
$$
where $\equiv$ denotes reduction by $f$.
Then for each $i < i_0$, the polynomial $g_{i+1}$ must be
added to the 
basis.
Computing the remaining $S$-polynomials,
$$
\eqalignno{
S(f,g_{i_0}) &= x^{9q-3} w^{3q-3} f - y z^2 g_{i_0} \cr
&= - x^{9q}yz w^{3q} + x^{6q-3}y^{q+5\over 2} z^4w^{7q-3 \over 2} \cr
&\equiv
- x^{6q-1}y^{q+3\over 2}z^2w^{7q-1 \over 2}
+ x^{6q-1}y^{q+3\over 2}z^2w^{7q-1 \over 2} = 0, \cr
}
$$
where $\equiv$ denotes reduction by $f$ and $g_{i_0}$,
and for $i<j$,
$$
\eqalignno{
S(g_i, g_j) &= x^{2(j-i)} w^{j-i} g_i - y^{j-i} g_j \cr
&= - x^{3q-2+2j}y^{2q+1-i}z^2w^{2q-1+j}
+ x^{3q-2+2j}y^{2q+1-i}z^2w^{2q-1+j} = 0. \cr
}
$$
Therefore
$\{f, g_0, \ldots, g_{i_0}\}$ is a
Gr\"obner basis of
$J + \Fq I$; this basis is also reduced.
Thus when $q \ge 3$,
the maximal $w$-degree of the Gr\"obner basis is $\xq={7q-3 \over 2}$,
the maximal total degree is $\tq=12q - 1$,
and the number of elements is $\cq={3q+3 \over 2}$.
\qed

In the following example we again use the ideals
$I$ and $J$ from the previous Proposition and
exchange their roles, in order to
exhibit periodic behavior of both
the cardinality function $\cq$ and the total degree function $\tq$
of the elements of the
reduced Gr\"obner basis of $I + \Fq J$,
with periodic behavior starting not with $q = 1$
but at the next level, at $q = p$.

\prop
\label{\lmbdperiod}
Let $R = \bbZ/3 \bbZ [x,y,z,w]$,
$I = (x^5y^2zw-xy^3z^2w)$, 
$J = (xy^2z^3w^2-x^3yzw^3)$ and $q=3^e$.
Using the
reverse lexicographic ordering
(with $w<z<x<y$)
the
reduced Gr\"obner basis of $I+J$ is
$$
\{x^5y^2zw-xy^3z^2w, xy^2z^3w^2-x^3yzw^3, x^7yzw^3-x^3y^2z^2w^3\};
$$
the
reduced Gr\"obner basis of $I + \Fq J$
for $q$ a positive even power of $3$ is
$$
\{x^5y^2zw-xy^3z^2w,
x y^{{9 \over 4}q-{1 \over 4}} z^{{13 \over 4}q-{1 \over 4}} w^{2q} -
x^3 y^{{7 \over 4}q-{3 \over 4}} z^{{7 \over 4}q-{3 \over 4}} w^{3q}\};
$$
and the
reduced Gr\"obner basis of $I + \Fq J$ for $q$ an odd power of $3$ is
$$
\eqalignno{
\{
x^5y^2zw-xy^3z^2w,
&\ x^3 y^{{9 \over 4}q-{3 \over 4}} z^{{13 \over 4}q-{3 \over 4}} w^{2q} -
x y^{{7 \over 4}q-{1 \over 4}} z^{{7 \over 4}q-{1 \over 4}} w^{3q}, \cr
& x y^{{9 \over 4}q+{1 \over 4}} z^{{13 \over 4}q+{1 \over 4}} w^{2q} -
x^3 y^{{7 \over 4}q-{1 \over 4}} z^{{7 \over 4}q-{1 \over 4}} w^{3q} \}.
}
$$
The corresponding functions are given by $\xq = 3q$, 
$$
\tq = \cases{12&if $q=1$, \cr
\noalign{\vskip2pt}
             {15 \over 2}q+{3 \over 2} &if $q=3^e$, $e$ odd, \cr
\noalign{\vskip4pt}
             {15 \over 2}q+{1 \over 2} &if $q=3^e$, $e>0$ even }
\quad {\rm and} \quad
\cq = \cases{3&if $q=1$ or $q=3^e$, $e$ odd, \cr
                2&if $q=3^e$, $e>0$ even }
$$
for all $q$.
\endb

\proof
As in the proof of Proposition~\lmbdperiodtwo, 
the case $q = 1$ can be computed directly by Macaulay2 and is
left to the reader.

Assume that $q \ge 3$; write $q=3^e$, with $e>0$.
Observe that whenever $a \ge 5$, $b \ge 2$, $c \ge 1$ and $d \ge 1$,
the monomial $x^a y^b z^c w^d$ reduces to
$x^{a-4} y^{b+1} z^{c+1} w^d$ modulo $I$,
so the normal form of
$x^a y^b z^c w^d$ modulo $I$
is $x^{a-4k} y^{b+k} z^{c+k} w^d$,
where $k$ is the largest integer such that
$a-4(k-1) \ge 5$, i.e.,
$a \ge 4k + 1$.
In reducing
the generator $x^qy^{2q}z^{3q}w^{2q} - x^{3q}y^qz^qw^{3q}$ of $\Fq J$
using the generator $g=x^5y^2zw-xy^3z^2w$ of $I$,
we need to consider the cases when $e$ is even and odd
separately.
If $e$ is even, then $q \equiv 1 (\mod 4)$, so
the monomial $x^qy^{2q}z^{3q}w^{2q}$ can be reduced
$k= {q-1 \over 4}$ times and
the monomial $x^{3q}y^qz^qw^{3q}$ can be reduced
$k={3q-3 \over 4}$ times using $g$.
Similarly, if $e$ is odd, then
$q \equiv 3 (\mod 4)$, so
the monomial $x^qy^{2q}z^{3q}w^{2q}$ can be reduced
$k= {q-3 \over 4}$ times and
the monomial $x^{3q}y^qz^qw^{3q}$ can be reduced
$k={3q-1 \over 4}$ times using $g$.
In particular,
when $q \ge 3$,
the generator $x^qy^{2q}z^{3q}w^{2q}-x^{3q}y^qz^qw^{3q}$ of $\Fq J$
reduces to $f'$,
where
$$
f' := \Bigl\lbrace
\matrix{x y^{{9 \over 4}q-{1 \over 4}} z^{{13 \over 4}q-{1 \over 4}} w^{2q} -
x^3 y^{{7 \over 4}q-{3 \over 4}} z^{{7 \over 4}q-{3 \over 4}} w^{3q},
& \hbox{\ $e$ even positive} \cr
\noalign{\vskip 3pt}
x^3 y^{{9 \over 4}q-{3 \over 4}} z^{{13 \over 4}q-{3 \over 4}} w^{2q} -
x y^{{7 \over 4}q-{1 \over 4}} z^{{7 \over 4}q-{1 \over 4}} w^{3q},
& \hbox{\ $e$ odd}. \hfill \cr
}
$$

If $e>0$ is even,
the S-polynomial of $f'$ and the generator $g$ of $I$ is
$$
\eqalignno{
S(g, f') &=
y^{{9 \over 4}q-{9 \over 4}} z^{{13 \over 4}q-{5 \over 4}} w^{2q-1} g
   - x^4 f'\cr
&= - x y^{{9 \over 4}q+{3 \over 4}} z^{{13 \over 4}q+{3 \over 4}} w^{2q}
+ x^7 y^{{7 \over 4}q-{3 \over 4}} z^{{7 \over 4}q-{3 \over 4}} w^{3q} \cr
&\equiv - x^3 y^{{7 \over 4}q+{1 \over 4}} z^{{7 \over 4}q+{1 \over 4}} w^{3q}
+ x^3 y^{{7 \over 4}q+{1 \over 4}} z^{{7 \over 4}q+{1 \over 4}} w^{3q} = 0,\cr
}
$$
where $\equiv$ denotes one reduction using $f'$ on the first term,
and one reduction using $g$ on the second term.
For the case in which
$e$ is positive and even, this proves that the basis in
the Proposition is a Gr\"obner basis; this basis is
reduced, and the three functions can be computed directly.

Finally assume that $q$ is an odd power of $p = 3$.
Then the S-polynomial of $f'$ and the generator $g$ of $I$ is
$$
\eqalignno{
S(g, f') &=
y^{{9 \over 4}q-{11 \over 4}} z^{{13 \over 4}q-{7 \over 4}} w^{2q-1} g -
x^2 f' =- x y^{{9 \over 4}q+{1 \over 4}} z^{{13 \over 4}q+{1 \over 4}} w^{2q}
+x^3 y^{{7 \over 4}q-{1 \over 4}} z^{{7 \over 4}q-{1 \over 4}} w^{3q}, \cr
}
$$
which is non-zero and reduced with respect to the set $\{f',g\}$.
So let
$$
h := x y^{{9 \over 4}q+{1 \over 4}} z^{{13 \over 4}q+{1 \over 4}} w^{2q} -
x^3 y^{{7 \over 4}q-{1 \over 4}} z^{{7 \over 4}q-{1 \over 4}} w^{3q}
$$
be added to $f'$ and $g$ in the procedure to form
a Gr\"obner basis of $I + \Fq J$.
Computing the remaining $S$-polynomials, we find that
$$
\eqalignno{
S(g,h) &= y^{{9 \over 4}q-{7 \over 4}} z^{{13 \over 4}q-{3 \over 4}} w^{2q-1}g
- x^4 h
= - x y^{{9 \over 4}q+{5 \over 4}} z^{{13 \over 4}q+{5 \over 4}} w^{2q}
  + x^7 y^{{7 \over 4}q-{1 \over 4}} z^{{7 \over 4}q-{1 \over 4}} w^{3q} \cr
&\equiv 
 - x^3 y^{{7 \over 4}q+{3 \over 4}} z^{{7 \over 4}q+{3 \over 4}} w^{3q} +
 x^3 y^{{7 \over 4}q+{3 \over 4}} z^{{7 \over 4}q+{3 \over 4}} w^{3q}= 0, \cr
}$$
where $\equiv$ denotes one reduction using $h$ on the first term
and one reduction using $g$ on the second term,
and the last $S$-polynomial
$$
\eqalignno{
S(f',h) &=y z f' - x^2 h
=- x y^{{7 \over 4}q+{3 \over 4}} z^{{7 \over 4}q +{3 \over 4}} w^{3q}
  + x^5 y^{{7 \over 4}q-{1 \over 4}} z^{{7 \over 4}q-{1 \over 4}} w^{3q} \cr
&\equiv -x y^{{7 \over 4}q+{3 \over 4}} z^{{7 \over 4}q +{3 \over 4}} w^{3q}
  + x y^{{7 \over 4}q+{3 \over 4}} z^{{7 \over 4}q +{3 \over 4}} w^{3q}
= 0, \cr
}
$$
where $\equiv$ denotes one reduction using $g$ on the second term.
Therefore $\{f', g, h\}$ is indeed a Gr\"obner basis of $I + \Fq J$
for $q$ an odd power of $p$.  This basis is also reduced, so
in this case the maximal $w$-degree of the Gr\"obner basis is $\xq=3q$,
the maximal total degree is
$\tq= {15 q + 3 \over 2}$,
and the number of elements is $\cq=3$.
\qed

In the next example we show that the function $\xq$ also 
can vary periodically.
In the example in Proposition~\lmbdperiod,
$\cq$ alternated between constant functions for the ideals $J + \Fq I$.
The next example shows that the function $\cq$ can
vary periodically between linear functions as well.
Moreover, the asymptotic patterns for all
three functions of the ideals $J + \Fq I$
begin further along, at $q=p^2$.

\prop
\label{\lmbdpudud}
Let $R = \bbZ/3\bbZ[x,y,z]$,
$I = (x^2y^2z-xyz^2)$,
$J = (xy^2z^5-x^2yz)$ and $q=3^e$.
Then with the reverse lexicographic ordering (with $z<y<x$) the
reduced Gr\"obner basis for
$J + I$ is
$$
\{xy^2z^5-x^2yz, x^2y^2z-xyz^2, xyz^6-x^3yz \};
$$
the
reduced Gr\"obner basis for
$J + I^{[3]}$ is
$$
\{xy^2z^5-x^2yz, x^6y^6z^3-x^4y^2z^2, x^7y^5z-x^4y^2z^4, 
x^8y^4z-x^5yz^4, x^5yz^8-x^9y^3z\};
$$
if $e \ge 2$ is even
the
reduced Gr\"obner basis for
$J + I^{[q]}$ is
$$
\eqalign{
\{&xy^2z^5-x^2yz,
x^{{9 \over 4}q-{1 \over 4}+k}y^{{7 \over 4}q+{1 \over 4}-k}z -
x^{{3 \over 2}q-{1 \over 2}+k}y^{{1 \over 2}q+{1 \over 2}-k}z^2, \cr
& x^{{11 \over 4}q+{1 \over 4}+j}y^{{5 \over 4}q-{1 \over 4}-j}z -
x^{2q-1}yz^{6+4j},
x^{2q-1}yz^{2q+4} -
x^{{13 \over 4}q-{1 \over 4}}y^{{3 \over 4}q+{1 \over 4}}z \cr
& ~|~ 0 \le k \le (q-1)/2, 0 \le j \le (q-3)/2\};
\cr}
$$
and if $e \ge 3$ is odd then
the
reduced Gr\"obner basis for
$J + \Fq I$ is
$$
\eqalign{
\{&xy^2z^5-x^2yz,
x^{{9 \over 4}q-{3 \over 4}}y^{{7 \over 4}q+{3 \over 4}}z^3 -
x^{{3 \over 2}q-{1 \over 2}}y^{{1 \over 2}q+{1 \over 2}}z^2,
 x^{{9 \over 4}q+{1 \over 4}+k}y^{{7 \over 4}q-{1 \over 4}-k}z -
x^{{3 \over 2}q-{1 \over 2}+k}y^{{1 \over 2}q+{1 \over 2}-k}z^4, \cr
& x^{{11 \over 4}q+{3 \over 4}+j}y^{{5 \over 4}q-{3 \over 4}-j}z -
x^{2q-1}yz^{8+4j},
x^{2q-1}yz^{2q+6} -
x^{{13 \over 4}q+{1 \over 4}}y^{{3 \over 4}q-{1 \over 4}}z \cr
& ~|~ 0 \le k \le (q-1)/2, 0 \le j \le (q-3)/2\}.
\cr}
$$
The associated functions are
$$
{
\xq = \cases{8&if $q=3$, \cr
                           2q+4&if $q=3^e$, $e  \ge 0$ even, \cr
                           2q+6&if $q=3^e$, $e \ge 3$ odd, \cr}
}
\qquad  \qquad
{
\tq = \cases{15&if $q=3$, \cr
                           4q+4&if $q=3^e$, $e  \ge 0$ even, \cr
                           4q+6&if $q=3^e$, $e \ge 3$ odd \cr}
}
$$
$$
{\rm and} \qquad \cq = \cases{5&if $q=3$,\cr
                          q+2&if $q=3^e$, $e  \ge 0$ even, \cr
                          q+3&if $q=3^e$, $e \ge 3$ odd \cr}
$$
for all $q$.

\endb

\proof
The Gr\"obner
bases for $J+I$ and $J+I^{[3]}$ can be computed with Macaulay2, and
are left to the reader.
For the rest of the proof, assume $q=p^e$ with $e \ge 2$.
Let $g=xy^2z^5-x^2yz$ be the generator of the ideal $J$.
We need to reduce the generator $x^{2q}y^{2q}z^q-x^qy^qz^{2q}$
of $\Fq I$ to normal form modulo $g$.
Observe that whenever $a \ge 1$, $b \ge 2$, and $c \ge 5$,
then $x^ay^bz^c$ reduces to $x^{a+1}y^{b-1}z^{c-4}$,
so the normal form of the monomial $x^ay^bz^c$ is
the monomial $x^{a+k}y^{b-k}z^{c-4k}$, where
$k$ is the largest integer such that
$b-(k-1) \ge 2$ and $c-4(k-1) \ge 5$; i.e.,
$b \ge k+1$ and $c \ge 4k+1$.
For the monomial $x^{2q}y^{2q}z^q$, 
$k$ is the largest integer such that
$2q \ge k+1$ and $q \ge 4k+1$; in this case,
if the latter inequality holds, then the former is true
as well, so we only need to find the largest
integer $k$ for which $q \ge 4k+1$.
If $e$ is even, then $q \equiv 1 \mod 4$,
so $k=(q-1)/4$, and the normal form of $x^{2q}y^{2q}z^q$
is
$x^{{9 \over 4}q-{1 \over 4}}y^{{7 \over 4}q+{1 \over 4}}z$.
If $e$ is odd, then
$q \equiv 3 \mod 4$, so $k=(q-3)/4$, and
the normal form of $x^{2q}y^{2q}z^q$
is
$x^{{9 \over 4}q-{3 \over 4}}y^{{7 \over 4}q+{3 \over 4}}z^3$.
Similarly, $x^qy^qz^{2q}$ reduces $k$ times using $g$ to
its normal form when $k$ is the largest integer such that
$q \ge k+1$ and $2q \ge 4k+1$.  As before we can ignore the
first inequality.  For all $e \ge 2$, we get $k=(2q-2)/4=(q-1)/2$,
so the normal form of $x^qy^qz^{2q}$ is
$x^{{3 \over 2}q-{1 \over 2}}y^{{1 \over 2}q+{1 \over 2}}z^2$.
The the normal form for the generator
$x^{2q}y^{2q}z^q-x^qy^qz^{2q}$
of $\Fq I$ is
$$
f' := \Bigl\lbrace
\matrix{x^{{9 \over 4}q-{1 \over 4}}y^{{7 \over 4}q+{1 \over 4}}z -
x^{{3 \over 2}q-{1 \over 2}}y^{{1 \over 2}q+{1 \over 2}}z^2,
& \hbox{\ $e$ even}, \cr
\noalign{\vskip 3pt}
x^{{9 \over 4}q-{3 \over 4}}y^{{7 \over 4}q+{3 \over 4}}z^3 -
x^{{3 \over 2}q-{1 \over 2}}y^{{1 \over 2}q+{1 \over 2}}z^2,
& \hbox{\ $e$ odd}. \hfill \cr
}
$$

Suppose that $e \ge 2$ is even.  Define the polynomials
$$
\eqalign{
f_k&:=x^{{9 \over 4}q-{1 \over 4}+k}y^{{7 \over 4}q+{1 \over 4}-k}z -
x^{{3 \over 2}q-{1 \over 2}+k}y^{{1 \over 2}q+{1 \over 2}-k}z^2 \qquad
{\rm for~} 0 \le k \le (q-1)/2, \cr
h_j&:=x^{{11 \over 4}q+{1 \over 4}+j}y^{{5 \over 4}q-{1 \over 4}-j}z -
x^{2q-1}yz^{6+4j} \qquad {\rm for~} 0 \le j \le (q-3)/2, {\rm ~and} \cr
r&:=x^{2q-1}yz^{2q+4} -
x^{{13 \over 4}q-{1 \over 4}}y^{{3 \over 4}q+{1 \over 4}}z. \cr
}$$
Note that $f'=f_0$.
When $ 0 \le k \le (q-3)/2$, the $S$-polynomial
$$
\eqalign{
S(g,f_k) &= x^{{9 \over 4}q-{5 \over 4}+k}y^{{7 \over 4}q-{7 \over 4}-k} g
  - z^4 f_k \cr
   &= -x^{{9 \over 4}q-{1 \over 4}+(k+1)}y^{{7 \over 4}q+{1 \over 4}-(k+1)}z
  + x^{{3 \over 2}q-{1 \over 2}+k}y^{{1 \over 2}q+{1 \over 2}-k}z^6 \cr
&\equiv -x^{{9 \over 4}q-{1 \over 4}+(k+1)}y^{{7 \over 4}q+{1 \over 4}-(k+1)}z
  + x^{{3 \over 2}q-{1 \over 2}+(k+1)}y^{{1 \over 2}q+{1 \over 2}-(k+1)}z^2
  = -f_{k+1}, \cr
}
$$
where $\equiv$ denotes a reduction using $g$ on the second term.
Therefore the polynomials $f_k$ for $ 0 \le k \le (q-1)/2$ are
included with $g$ and $f'$ in the procedure to compute the
Gr\"obner basis.  The last polynomial in this family is
$f_{(q-1)/2}=x^{{11 \over 4}q-{3 \over 4}}y^{{5 \over 4}q+{3 \over 4}}z -
x^{2q-1}yz^{2}$.  Then
$$
\eqalign{
S(g,f_{(q-1)/2})
  & = x^{{11 \over 4}q-{7 \over 4}}y^{{5 \over 4}q-{5 \over 4}} g
   - z^4f_{(q-1)/2} \cr
  &= -x^{{11 \over 4}q+{1 \over 4}}y^{{5 \over 4}q-{1 \over 4}}z
     + x^{2q-1}yz^{6} = -h_0. \cr
} $$
Similarly, the $S$-polynomial $S(g,h_j)=-h_{j+1}$ for all
$ 0 \le j \le (q-5)/2$, so the polynomials $h_j$ for
$ 0 \le j \le (q-3)/2$ are appended to the basis.  The final
polynomial in this list is
$h_{(q-3)/2} = x^{{13 \over 4}q-{5 \over 4}}y^{{3 \over 4}q+{5 \over 4}}z
  - x^{2q-1}yz^{2q}$.  Then
$$
\eqalign{
S(g,h_{(q-3)/2}) 
  &= x^{{13 \over 4}q-{9 \over 4}}y^{{3 \over 4}q-{3 \over 4}} g
   - z^4 h_{(q-3)/2}
= -x^{{13 \over 4}q-{1 \over 4}}y^{{3 \over 4}q+{1 \over 4}}z
     + x^{2q-1}yz^{2q+4} = r. \cr
}$$
Therefore $r$ is also added to the basis by the Buchberger
algorithm.  All of the
remaining $S$-polynomials reduce to zero modulo this
set of polynomials, 
so the set
$\{g,f_k,h_j,r ~|~ 0 \le k \le (q-1)/2, 0 \le j \le (q-3)/2\}$ is
a Gr\"obner basis for $J + \Fq I$ in the case that $e \ge 2$ is even.

Finally, suppose that $e \ge 3$ is odd.
We have already shown that the polynomials
$g =xy^2z^5-x^2yz$ and
$f'=x^{{9 \over 4}q-{3 \over 4}}y^{{7 \over 4}q+{3 \over 4}}z^3 -
x^{{3 \over 2}q-{1 \over 2}}y^{{1 \over 2}q+{1 \over 2}}z^2$
are a basis for $J+\Fq I$.
Define the polynomials
$$
\eqalign{
s_k&:= x^{{9 \over 4}q+{1 \over 4}+k}y^{{7 \over 4}q-{1 \over 4}-k}z -
x^{{3 \over 2}q-{1 \over 2}+k}y^{{1 \over 2}q+{1 \over 2}-k}z^4
\qquad {\rm for~} 0 \le k \le (q-1)/2, \cr
t_j&:=x^{{11 \over 4}q+{3 \over 4}+j}y^{{5 \over 4}q-{3 \over 4}-j}z -
x^{2q-1}yz^{8+4j} \qquad {\rm for~} 0 \le j \le (q-3)/2, {\rm ~and} \cr
u&:=x^{2q-1}yz^{2q+6} -
x^{{13 \over 4}q+{1 \over 4}}y^{{3 \over 4}q-{1 \over 4}}z. \cr
}$$
By an argument very similar to the proof above, we
get that
$S(g,f')=-s_0$ and
$S(g,s_k) \equiv -s_{k+1}$ for all $ 0 \le k \le (q-3)/2$, where
$\equiv$ denotes a reduction by $g$.  Then
$S(g,s_{(q-1)/2}) = t_0$, and
$S(g,t_j) = -t_{j+1}$ when $0 \le j \le (q-5)/2$.  Taking
one further $S$-polynomial with $g$,
$S(g,t_{(q-3)/2}) = u$.   Finally, all of the remaining
$S$ polynomials reduce to 0 modulo these polynomials, 
so
the set $\{g,f',s_k,t_j,u ~|~ 0 \le k \le (q-1)/2, 0 \le j \le (q-3)/2\}$
is a Gr\"obner basis for $J + \Fq I$ when $e \ge 3$ is odd.

Since in each case the Gr\"obner basis we computed is also reduced,
the results on the functions associated to these ideals then
follow immediately from these bases.
\qed

If we change the characteristic in Proposition \lmbdpudud\
to $p=2$, we find that 
the $x_n$-degree function $\xq$ is dependent on
the characteristic of the field~$F$ as well; in fact,
all three functions $\xq$, $\tq$ and
$\cq$ are altered, and the periodicity is lost.

\prop
\label{\lmbdchar}
Let $R = \bbZ/2\bbZ[x,y,z]$,
$I = (x^2y^2z-xyz^2)$, 
$J = (xy^2z^5-x^2yz)$ and
$q=2^e$.
Then with the reverse lexicographic ordering (with $z<y<x$) the
reduced Gr\"obner basis for
$J+\Fq I$ with $e \ge 3$ is
$$
\eqalign{
\{&xy^2z^5-x^2yz,
x^{{9 \over 4}q-1}y^{{7 \over 4}q+1}z^4 -
x^{{3 \over 2}q-1}y^{{1 \over 2}q+1}z^4,
x^{{9 \over 4}q+j}y^{{7 \over 4}q-j}z -
x^{{3 \over 2}q+j}y^{{1 \over 2}q-j}z, \cr
& x^{{11 \over 4}q+k}y^{{5 \over 4}q-k}z -
x^{2q-1}yz^{5+4k},
x^{2q-1}yz^{2q+5} - x^{{13 \over 4}q-1}y^{{3 \over 4}q+1}z^5 \cr
& ~|~ 0 \le j \le (q-2)/2, 0 \le k \le (q-2)/2\}.
\cr}
$$
The associated functions satisfy
$\xq=2q+5$, $\tq=4q+5$ and $\cq = q+3$ for
$q \ge 2^3$.
\endb

\proof
Assume that $q=2^e$ with $e \ge 3$.
The main difference between the proof of this Proposition and that
of Proposition \lmbdpudud\ lies in the reduction of the generator
$x^{2q}y^{2q}z^q-x^qy^qz^{2q}$ of $\Fq I$ modulo the
generator $g=xy^2z^5-x^2yz$ of $J$.
As in the earlier proof, if $a \ge 1$, $b \ge 2$, and $c \ge 5$, then
the normal form of the monomial $x^ay^bz^c$ is
the monomial $x^{a+k}y^{b-k}z^{c-4k}$, where
$k$ is the largest integer such that
$b \ge k+1$ and $c \ge 4k+1$.
For the monomial $x^{2q}y^{2q}z^q$, we get that
$k$ is the largest integer such that
$q \ge 4k+1$, i.e. $k \le (q-1)/4$. Since $q$ is even,
and moreover divisible by 4, this means that $k=(q-4)/4$.
Then the normal form of $x^{2q}y^{2q}z^q$ is
$x^{{9 \over 4}q-1}y^{{7 \over 4}q+1}z^4$.  Similarly,
to find the normal form of $x^qy^qz^{2q}$ we reduce
the monomial $k$ times where $k$ is the largest
integer satisfying $2q \ge 4k+1$, or $k \le (2q-1)/4$.
Then $k=(2q-4)/4$.  So the generator of $\Fq I$ reduces
to the polynomial
$f':=x^{{9 \over 4}q-1}y^{{7 \over 4}q+1}z^4 -
x^{{3 \over 2}q-1}y^{{1 \over 2}q+1}z^4$ modulo $g$.

Next define the polynomials
$f_j:=x^{{9 \over 4}q+j}y^{{7 \over 4}q-j}z -
x^{{3 \over 2}q+j}y^{{1 \over 2}q-j}z$
for $0 \le j \le (q-2)/2$,
$h_k:=x^{{11 \over 4}q+k}y^{{5 \over 4}q-k}z -
x^{2q-1}yz^{5+4k}$ for $0 \le k \le (q-2)/2$,
and
$r=:x^{2q-1}yz^{2q+5} - x^{{13 \over 4}q-1}y^{{3 \over 4}q+1}z^5$.
Following steps very similar to those in the proof
of Proposition \lmbdpudud, we find that
$S(g,f') \equiv -f_0$ and
$S(g,f_j) \equiv -f_{j+1}$ for $0 \le j \le (q-4)/2$,
where $\equiv$ denotes a single reduction by $g$ in each case.
Also,
$S(g,f_{(q-2)/2})=-h_0$, $S(g,h_k)=-h_{k+1}$
for $0 \le k \le (q-4)/2$, and $S(g,h_{(q-2)/2})=r$.
All other $S$-polynomials reduce to zero with this
basis, 
and no element of this set can be reduced by any other,
so the set $\{g,f',f_j,h_k,r ~|~ 0 \le j \le (q-2)/2,
0 \le k \le (q-2)/2\}$
is a reduced Gr\"obner basis.
\qed

In the final example we show that
it need not be the case that the total degree of the Gr\"obner basis
of $J + \Fq I$ is bounded above
by $q \cdot \max\{{\rm Gbdeg~} I, {\rm Gbdeg~} J\}$, where
${\rm Gbdeg}$ denotes the total degree of the 
reduced Gr\"obner basis
(with the reverse lexicographic ordering).

\prop
\label{\lmbdcoef}
Let $R = \bbZ/3\bbZ[x,y,z,w]$.
The ideal $J + \Fq I$ with
$I = (x^2y^2zw^5-xyz^2w^2)$, $J = (xy^2z^3w-xyzw^3)$ and $q=3^e$
has the reduced Gr\"obner basis
$$
\lbrace
xy^2z^3w-xyzw^3,
 x^{2q}y^{(3q+1-2k)/2}zw^{6q-1+2k}-x^qyz^{2k+2}w^{4q-2}
~|~ 0 \le k \le (3q-1)/2
\rbrace
$$
with respect to the reverse lexicographic ordering with $w<z<y<x$.
Therefore
the maximal $w$-degree of the Gr\"obner basis is $\xq=9q-2$,
the maximal total degree is $\tq=11q$, and 
the number of elements is $\cq=3(q+1)/2$ for all $q$.
\endb

Therefore $q \cdot \max\{{\rm Gbdeg~} I, {\rm Gbdeg~} J\} = 
q \cdot \max\{10, 7\} < 11q
= {\rm Gbdeg~} (J + \Fq I)$.

\medskip

\proof
Define $f:=xy^2z^3w-xyzw^3$ and for $0 \le k \le (3q-1)/2$, define
$$
g_k:=x^{2q}y^{(3q+1-2k)/2}zw^{6q-1+2k}-x^qyz^{2k+2}w^{4q-2}.
$$

Note that in this example the generator
$x^{2q}y^{2q}z^{q}w^{5q}-x^{q}y^{q}z^{2q}w^{2q}$
of $\Fq I$ is not in normal form modulo $J$.  The term
$x^{2q}y^{2q}z^{q}w^{5q}$ reduces (after $(q-1)/2$ reductions)
to $x^{2q}y^{(3q+1)/2}zw^{6q-1}$, and
$x^{q}y^{q}z^{2q}w^{2q}$ reduces (using $q-1$ reductions)
to $x^qyz^{2}w^{4q-2}$, resulting in the 
basis element
$
g_0=x^{2q}y^{(3q+1)/2}zw^{6q-1}-x^qyz^{2}w^{4q-2}.
$
The $S$-polynomial $S(f,g_0)$ is
$$
S(f,g_0)=-x^{2q}y^{(3q+1-2)/2}zw^{6q-1+2}+x^qyz^{2+2}w^{4q-2}=-g_1.
$$
Repeating this, when $k<(3q-1)/2$,
the $S$-polynomial $S(f,g_k)$ equals $-g_{k+1}$.
The remaining $S$-polynomials $S(f,g_{(3q-1)/2})$ and
$S(g_j,g_k)$ all
reduce to 0.
Thus the polynomials $f$ and $g_k$ for 
$0 \le k \le (3q-1)/2$ form a Gr\"obner basis.
This Gr\"obner basis is reduced, giving the
results on the functions.
\qed

\vskip 2ex
\vfill\eject

\centerline{\bf Summary table}
\vskip 4ex

\halign{\mvrule \vadjust{\vskip-2pt}\ #\hfill \mvrule
         & \ # \hfill \mvrule & \ # \hfill \mvrule
       & \ # \hfill \mvrule\cr
\noalign{\hrule}
{\bf Example}   & {\bf $\xq$}
        & {\bf $\tq$} & $\cq$ \cr
  &
        & &  \cr
\noalign{\vskip 2pt\hrule\vskip-2pt}
Prop.~\lmbdunb
       & linear & linear & constant \cr
\noalign{\vskip 2pt\hrule\vskip-2pt}
Prop.~\lmbdunbtwo     & linear & linear & linear \cr
\noalign{\vskip 2pt\hrule\vskip-2pt}
Prop.~\lmbdbd(a)
       & constant & linear & constant \cr
\noalign{\vskip 2pt\hrule\vskip-2pt}
Prop.~\lmbdbd(b) & linear & linear ($q \ge p$) & constant ($q \ge p$)\cr
\noalign{\vskip 2pt\hrule\vskip-2pt}
Prop.~\lmbdperiodtwo       & linear  ($q \ge p$)&
            linear ($q \ge p$) & linear \cr
\noalign{\vskip 2pt\hrule\vskip-2pt}
Prop.~\lmbdperiod
       & linear &
                 periodically & periodically  \cr
   & & {\ \ \ \ linear ($q \ge p$)} & {\ \ \ \ constant ($q \ge p$)}\cr
\noalign{\vskip 2pt\hrule\vskip-2pt}
Prop.~\lmbdpudud
       & periodically &
                 periodically & periodically \cr
   & \ \ \ \ linear ($q \ge p^2$)  & \ \ \ \ linear ($q \ge p^2$)
          & \ \ \ \ linear ($q \ge p^2$)  \cr
\noalign{\vskip 2pt\hrule\vskip-2pt}
Prop.~\lmbdcoef
       & linear &
                 linear (high coeff.) & linear  \cr
\noalign{\vskip 2pt\hrule\vskip-2pt}
}

\vskip 4ex

All of these examples satisfy Katzman's conjecture
that the $x_n$-degree $\xq$ of the
reduced Gr\"obner basis
of $J + \Fq I$ is bounded above linearly in $q$.
Furthermore,
in all of these examples the total degree and cardinality
of the Gr\"obner basis
are also bounded above linearly in $q$.
However, we are left with the open question of whether 
the behavior of the functions $\xq$, $\tq$ and $\cq$ 
(eventually) follows one of the patterns in the table above,
and 
whether linear
upper bounds on $\xq$, $\tq$ and $\cq$ 
hold, for all 
ideals $I$ and $J$ in a polynomial ring.

\vskip 0.3truecm\relax

\leftline{\bf Acknowledgment}

\vskip 0.2truecm\relax

The authors thank Aldo Conca and Enrico Sbarra for
helpful discussions on an earlier version of this paper.

\vskip 0.3truecm\relax

\leftline{\bf Appendix: Macaulay2 code}

\vskip 0.2truecm\relax

We used variations of the following Macaulay2 code for our calculations,
included for the readers interested in making further computations.

{\bigskip
\bgroup
\parindent=1ex
\narrower
\tt
\obeylines
\obeyspaces
\baselineskip=0pt
{\bf Input: polynomial ring $R$, ideals $I$, $J$}
{\bf Output: fn(e) = Gr\"obner basis of $\bf J + I^{[p^e]}$,}
{\bf \ \hphantom{Output:} df(e) = maximal total degree of an element of the Gr\"obner basis.}
p = 3
R = ZZ/p[x,y,z,MonomialSize=>16];
I = ideal(y\^ 2*z-x\^ 2);
J = ideal(y\^ 3-x*y);
fn = e -> (transpose gens gb (J+I\^ (p\^ e)))
df = e -> (L = $\{ \}$; i = 0;
\           G = gens gb (J + I\^ (p\^ e));
\           l = rank source G;
\           while i < l do (
\                 L = prepend (degree G\_(0,i), L);
\                 i = i + 1; );
\           max L)
\ %

\egroup
}

\vskip 0.3truecm\relax

\leftline{\bf References}

\vskip 0.2truecm\relax

\font\eightrm=cmr8 \def\rm{\fam0\eightrm}
\font\eightit=cmti8 \def\it{\fam\itfam\eightit}
\font\eightbf=cmbx8 \def\bf{\fam\bffam\eightbf}
\font\eighttt=cmtt8 \def\tt{\fam\ttfam\eighttt}
\rm
\baselineskip=9.9pt
\parindent=3.6em
\bgroup \catcode`\.=11

\item{[BS1]}
D. Bayer and M. Stillman,
A criterion for detecting $m$-regularity,
{\it Invent. Math.} {\bf 87} (1987), 1-11.

\item{[BS2]}
D. Bayer and M. Stillman,
A theorem on refining division orders by the reverse lexicographic order,
{\it Duke J. Math.} {\bf 55} (1987), 321-328.

\item{[CLO]}
D. Cox, J. Little and D. O'Shea,
{\it Ideals, Varieties, and Algorithms:
An introduction to computational algebraic geometry and commutative algebra},
Undergraduate Texts in Mathematics, Springer-Verlag, New York, 1992.

\item{[E]}
D. Eisenbud,
{\it Commutative Algebra with a View toward Algebraic Geometry},
Graduate Texts in Mathematics 150, Springer-Verlag, New York, 1995.

\item{[GS]}
D. R. Grayson and M. E. Stillman,
Macaulay2, a software system for research in algebraic geometry,
available at http://www.math.uiuc.edu/Macaulay2.

\item{[HH]}
M. Hochster and C. Huneke,
Tight closure, invariant theory, and the Brian\c con-Skoda Theorem,
{\it J. Amer. Math. Soc.} {\bf 3} (1990), {31-116}.

\item{[K]}
M. Katzman,
The complexity of Frobenius powers of ideals,
{\it J. Algebra} {\bf 203} (1998), 211-225.

\item{[P]}
K. Pardue,
Deformation classes of graded modules and maximal Betti numbers,
{\it Illinois J. Math} {\bf 40} (1996), 564-585.

\item{[S]}
K. E. Smith,
Tight closure commutes with localization in binomial rings,
{\it Proc. Amer. Math. Soc.} {\bf 129} (2001), 667-669.

\item{[Sw]}
I. Swanson,
Powers of ideals: Primary decompositions, Artin-Rees lemma and regularity,
{\it Math. Ann.} {\bf 307} (1997), 299-313.

\egroup

\bigskip

\noindent Department of Mathematics and Statistics,
University of Nebraska, Lincoln, Nebraska 68588-0323, USA,
E-mail: {\tt smh@math.unl.edu}

\medskip

\noindent Department of Mathematical Sciences, New Mexico State University,
Las Cruces, New Mexico 88003-8001, USA,
E-mail: {\tt iswanson@nmsu.edu}

\end